\documentclass{amsart}
\usepackage{epsfig}
\usepackage{graphicx}
\usepackage{amsfonts}
\usepackage{amssymb}
\usepackage{qtree}
\usepackage{amsmath,latexsym,amsthm}

\newtheorem{theorem}{Theorem}[section]
\newtheorem{lemma}[theorem]{Lemma}

\newtheorem{corollary}[theorem]{Corollary}

\theoremstyle{definition}
\newtheorem{definition}[theorem]{Definition}
\newtheorem{example}[theorem]{Example}

\theoremstyle{remark}

\begin{document}

\title[$2$-adic valuation]{A binary tree representation for the 
$2$-adic valuation of a 
sequence arising from a rational integral}

\author{Xinyu Sun}
\address{Department of Mathematics,
Tulane University, New Orleans, LA 70118}
\email{xsun1@tulane.edu}

\author{Victor H. Moll}
\address{Department of Mathematics,
Tulane University, New Orleans, LA 70118}
\email{vhm@math.tulane.edu}


\date{\today}


\begin{abstract}
We analyze properties of the $2$-adic valuation of an integer
sequence that originates from an explicit 
evaluation of a quartic integral. We present a tree that encodes
this valuation. 
\end{abstract}

\maketitle

\newcommand{\realpart}{\mathop{\rm Re}\nolimits}
\newcommand{\imagpart}{\mathop{\rm Im}\nolimits}

\numberwithin{equation}{section}

\section{Introduction} \label{intro} 
\setcounter{equation}{0}

The integral
\begin{equation}
N_{0,4}(a;m) = \int_{0}^{\infty} \frac{dx}{(x^{4}+2ax^{2} + 1)^{m+1}},
\label{nzero4}
\end{equation}
with $a > -1$ is given by 
\begin{equation}
N_{0,4}(a,m) = \frac{\pi}{2} \frac{P_{m}(a)}{[2(a+1)]^{m+ 1/2}}
\label{relation}
\end{equation}
\noindent
where
\begin{equation}
P_{m}(a) = \sum_{l=0}^{m} d_{l,m} a^{l}
\end{equation}
\noindent
with 
\begin{equation}
d_{l,m} = 2^{-2m} \sum_{k=l}^{m} 2^{k} \binom{2m-2k}{m-k} \binom{m+k}{m}
\binom{k}{l}, \, \quad 0 \leq l \leq m. 
\end{equation}
\noindent
The reader will find in \cite{amram} a survey of the different proofs of 
(\ref{relation}) and an introduction to the many issues 
involved in the evaluation 
of definite integrals in \cite{moll-notices}.  

The study of combinatorial aspects of the sequence $d_{l}(m)$ was initiated 
in \cite{bomouni1} where the authors show that $d_{l}(m)$ form a {\em unimodal}
sequence, that is,
there exists and index $l^{*}$ such that 
$d_{0,m} \leq \ldots \leq d_{l^{*},m}$ and 
$d_{l^{*},m} \geq \ldots \geq d_{m,m}$. The fact that $d_{l,m}$  satisfies 
the stronger condition of {\em logconcavity} $d_{l-1,m}d_{l+1,m} \leq 
d_{l,m}^{2}$ has been recently established in \cite{kauers-paule}.

We consider here 
arithmetical properties of the sequence $d_{l,m}$. It is more convenient to 
analyze the auxiliary sequence 
\begin{equation}
A_{l,m} = l! \, m! \, 2^{m+l} d_{l,m}
 =  \frac{l! \, m!}{2^{m-l}} 
\sum_{k=l}^{m} 2^{k} \binom{2m-2k}{m-k} \binom{m+k}{m} 
\binom{k}{l} 
\label{intseq}
\end{equation}
\noindent
for $m \in \mathbb{N}$ and $0 \leq l \leq m$. The integral (\ref{nzero4}) 
is then given explicitly as 
\begin{equation}
N_{0,4}(a;m) = \frac{\pi}{\sqrt{2} \, m! \, (4(2a+1))^{m+1/2}} 
\sum_{l=0}^{m} A_{l,m} \frac{a^{l}}{l!}. 
\end{equation}

We present here a binary tree that encodes the 
$2$-adic valuation of 
$A_{l,m}$. Recall that, for $x \in \mathbb{N}$, the 
$2$-adic valuation $\nu_{2}(x)$ 
is the highest power of $2$ that divides $x$. This is extended to
$x  = a/b \in 
\mathbb{Q}$ via $\nu_{2}(x) = \nu_{2}(a) - \nu_{2}(b)$.  \\

Given $l \in \mathbb{N}$ we associate a tree $T(l)$, the {\em decision tree of } $l$, that provides a combinatorial interpretation of $\nu_{2}(A_{l,m})$. It 
has the following properties: \\

\noindent
1) Aside from the labels on the vertices, $T(l)$ depends only on the 
odd part of $l$. Therefore 
it suffices to consider $l$ odd. \\

\noindent
2) For $l$ odd, define $k^{*}(l) = \lfloor{ \log_{2} l \rfloor}$. The index 
$k^{*}$ is determined by $2^{k^{*}} < l < 2^{k^{*} + 1}$. \\

\noindent
3) The generations are labelled starting at $0$; that is, the root is generation
$0$. For $0 \leq k \leq k^{*}$, the $k$-th generation consists  
of $2^{k}$ vertices. These form a 
complete binary tree. \\

\noindent
4) A vertex with degree $1$ is called {\em terminal}. The edge containing a 
terminal vertex is called a {\em terminal branch}. The $k^{*}$-th generation 
contains $2^{k^{*}+1}-l$ terminal vertices. The tree $T(l)$ has one more 
generation consisting of $2(l-2^{k^{*}})$ terminal vertices. \\

\noindent
5) Each terminal vertex of $T(l)$ has a {\em vertex constant} attached to it.
These are given in Lemmas \ref{lem_1st} and \ref{lem_2nd}. Each
non-terminal vertex has two {\em edge functions} attached to it. \\

The main results presented here is: \\

\begin{theorem}
Let $l \in \mathbb{N}$. The data described above provides an explicit 
formula for the $2$-adic valuation of the sequence $A_{l,m}$.
\end{theorem}

The complete results are 
described in Section \ref{sec-tree} and illustrated here for $l=3$. \\

The sequence $\{ \nu_{2}(A_{3,m}): m \geq 1 \}$ satisfies 
$\nu_{2}(A_{3,2m-1}) = \nu_{2}(A_{3,2m})$. Therefore the 
subsequence $A_{3,2m+1}$, denoted by $C_{3,m}$, contains 
all the $2$-adic information 
of $A_{3,m}$. 
{{
\begin{figure}[ht]
\begin{center}
\centerline{\psfig{file=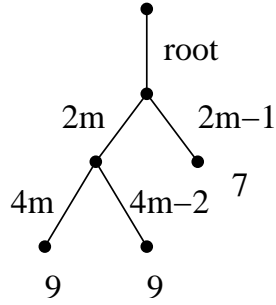,width=10em,angle=0}}
\caption{The decision tree for $l=3$ }
\label{dec-3}
\end{center}
\end{figure}
}}

For instance, at the first level in Figure \ref{dec-3} we have
two edges with functions $2m$ and $2m-1$ and the first generation consists
of two vertices, one of which is 
terminal with vertex constant $7$. This vertex is adjacent to the terminal 
branch labeled $2m-1$.  \\

This tree produces a formula for $\nu_{2}(C_{3,m})$ by the following 
mechanism: define
\begin{equation}
f_{3}(m) = \begin{cases}
7 +  \nu_{2} \left( \frac{m+1}{2} \right)  & \text{ if } m \equiv 1 \bmod 2, \\
9 +  \nu_{2} \left( \frac{m}{4} \right)  & \text{ if } m \equiv 0 \bmod 4, \\
9 +  \nu_{2} \left( \frac{m+2}{4} \right)  & \text{ if } m \equiv 2 \bmod 4. \\
\end{cases}
\end{equation}
\noindent
There is one expression per terminal branch. The 
numbers $9, \, 9, \, 7$ are the vertex constants of $T(3)$ 
and the arguments of $\nu_{2}$ in 
$f$ come from the branch labels. The tree now encodes the formula
\begin{equation}
\nu_{2}( C_{3,m} ) = f_{3}(m), \text{ for } m \geq 1.
\end{equation}

\section{The tree} \label{sec-tree} 
\setcounter{equation}{0}

In this section we describe a binary tree that encodes the $2$-adic valuation 
of the sequence $A_{l.m}$. This value is 
linked to that of the Pochhammer symbol
\begin{equation}
(a)_{n} := 
\begin{cases}
a(a+1)(a+2) \cdots (a+n-1), & \text{ for } n > 0 \\
1, & \text{ for } n =0, 
\end{cases}
\end{equation}
\noindent
via the identity 
\begin{equation}
\nu_{2}(A_{l,m}) = \nu_{2}((m+1-l)_{2l}) + l,
\label{poch}
\end{equation}
\noindent
established in \cite{amm1}. This is a generalization of the main 
result in \cite{bomosha}, namely,
\begin{equation}
\nu_{2}(A_{1,m}) = \nu_{2}(m(m+1)) + 1.
\end{equation}
\noindent
The expression (\ref{poch}) can also be written as 
\begin{equation}
\nu_{2}(A_{l,m}) = l + \sum_{j=-l+1}^{l} \nu_{2}(m+j). 
\label{poch-1}
\end{equation}

\medskip

To encode the information about $\nu_{2}(A_{l,m})$ we employ the notion of 
{\em simple sequences}. 

\begin{definition}
A sequence $\{ a_{n} : \, n \in \mathbb{N} \}$ is  called $s$-simple if 
there exists a number $s$ such that, for each 
$t \in \{ 0, \, 1, \, 2, \, \cdots \}$, we have
\begin{equation}
a_{st+1} = a_{st+2}  = \cdots = a_{s(t+1)}. 
\end{equation}
\end{definition}

In pictorial terms, $s$-simple 
sequences are formed by blocks of 
length $s$ where they attain the same value. In \cite{amm1} it is shown 
that, for fixed $l \in \mathbb{N}$,  the 
sequence $\{ \nu_{2}(A_{l,m}): m \geq l \}$ is 
$2^{1+ \nu_{2}(l)}-$simple. For instance,
\begin{equation}
\nu_{2}(A_{2,m}) = \{ 5, \, 5,\, 5,\, 5,\, 6, \, 6, \, 6, \, 6, \, 5, \, 
5, \, 5, \, 5, \, 7, \, 7, \, 7, \, 7, \, 5, \, 5, \, 5, \, 5, \ldots \}.
\end{equation}
\noindent
is $4$-simple.

\begin{definition}
Let $l \in \mathbb{N}$ be fixed. Define 
\begin{equation}
C_{l,m} = A_{l,l+(m-1) \cdot 2^{1+ \nu_{2}(l)}},
\end{equation}
\noindent
so that the  sequence
$\{ C_{l,m}: m \geq 1 \}$ reduces each block of $A_{l,m}$ to a 
single point. In particular, for $l$ odd we have $C_{l,m} = A_{l,2(m-1)}$.
\end{definition}

\noindent
{\bf The tree associated to} $l$.  We associate to each 
index $l \in \mathbb{N}$ a tree by the following rule: start with a 
{\em root} vertex. This root is the $0$-th generation of $T(l)$.
To the root vertex we attach the sequence 
\begin{equation}
\{ \nu_{2}( C_{l,m} ): \, m \geq 1 \}
\end{equation}
and ask whether 
\begin{equation}
\nu_{2}(C_{l,m}) - \nu_{2}(m) 
\end{equation}
is independent of $m$. If the answer is yes, we label the vertex $v_{0}$ 
with this constant value. This is the case for $l=4$ as shown in 
Figure \ref{dec-4}. If the 
answer is negative, we split the integers into classes modulo $2$
and create a vertex for each class. These two classes are attached to
two new vertices 
\begin{equation}
v_{1} \mapsto \{ \nu_{2}( C_{l,2m-1} ): \, m \geq 1 \}
\nonumber
\end{equation}
\noindent
and 
\begin{equation}
v_{2} \mapsto \{ \nu_{2}(C_{l,2m} ): \, m \geq 1 \}. 
\nonumber
\end{equation}
{{
\begin{figure}[ht]
\begin{center}
\centerline{\psfig{file=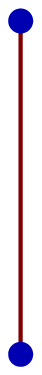,width=15em,angle=0}}
\caption{The decision tree for $l=4$ }
\label{dec-4}
\end{center}
\end{figure}
}}
\noindent 
Each positive answer produces the end of the branch and each negative one 
yields two new branches that need to be tested. The process stops 
when there are no more vertices that need to be tested.   \\

\noindent
{\bf Note}. Assume the vertex $v$ corresponding to the sequence 
$\{ 2^{k}(m-1) + a: \, m \geq 1 \}$ produces a negative answer. Then it 
splits in the next generation into two vertices corresponding to the 
sequences $\{ 2^{k+1}(m-1) + a: \, m \geq 1 \}$ and 
$\{ 2^{k+1}(m-1) + 2^{k} + a: \, m \geq 1 \}$. For 
instance, in Figure \ref{tree5}, the
vertex corresponding to $\{ 4m: \, m \geq 1 \}$, that is not terminal, splits 
into $\{ 8m: \, m \geq 1 \}$ and 
$\{ 8m + 4: \, m \geq 1 \}$. These two edges lead to terminal vertices. 
Theorem \ref{thm_main} shows that this 
process ends in a finite number of steps.  \\

{{
\begin{figure}[ht]
\begin{center}
{\centerline{\psfig{file=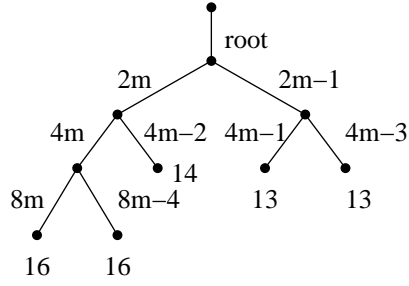,width=15em,angle=0} }}
\caption{The decision tree for $l=5$}
\label{tree5}
\end{center}
\end{figure}
}}

Figure \ref{dec-4} shows the decision tree for $l=4$ and Figure \ref{dec-3}
the corresponding one for $l=3$. In the latter figure the branches are labelled 
according to the arithmetic sequence they represent. 

At the second level we find the first appearance of a terminal vertex, namely
the one corresponding the edge marked $2m-1$. Its vertex constant is $7$, 
stating that 
\begin{equation}
\nu_{2}( C_{3,2m-1}) = \nu_{2}(m) + 7.
\end{equation}

\medskip

The first step in the analysis of the tree $T(l)$ is to reduce it to the 
case where $l$ is odd. 

\begin{theorem}
\label{odd-red}
The tree of an integer depends only upon its odd part, that is, for $b$ odd and 
$a \in \mathbb{N}$, we have
\begin{equation}
T(b) = T(2^{a} \, b),
\end{equation}
\noindent
up to relabel.
\end{theorem}

The proof of this theorem is based on a relation of the $2$-adic valuations 
of $C_{2l,m}$ and $C_{l,m}$. We establish first an auxiliary result for 
$A_{l,m}$. 

\begin{lemma}
Let $l, \, m \in \mathbb{N}$. Then 
\begin{equation}
\nu_{2}(A_{2l,2m}) = \nu_{2}(A_{l,m}) + 3l. 
\label{even-odd}
\end{equation}
\end{lemma}
\begin{proof}
The result is equivalent to 
\begin{equation}
\nu_{2} \left( a_{l}/a_{-l} \right) = 2l,
\label{new1}
\end{equation}
\noindent
where
\begin{equation}
a_{k} = \frac{ (2m+2k)!}{(m+k)!}. 
\end{equation}
\noindent 
This follows from 
\begin{equation}
\nu_{2} \left( a_{l}/a_{-l} \right) = 
\sum_{k=1-l}^{l} \nu_{2} \left( a_{k}/a_{k-1} \right),
\label{sum1}
\end{equation}
\noindent
and $a_{k}/a_{k-1} = 2(2m+2k-1)$, so that each term in the sum (\ref{sum1})
is equal to $1$.
\end{proof}

\begin{corollary}
\label{c-even}
Let $l, \, m \in \mathbb{N}$. Then 
\begin{equation}
\nu_{2} \left( C_{2l,m} \right) =  
\nu_{2} \left( C_{l,m} \right)  + 3l.
\end{equation}
\end{corollary}
\begin{proof}
The result follows from the identity
\begin{equation}
C_{2l,m} = A_{2l,2 \left[ l + (m-1) \cdot 2^{1+ \nu_{2}(l) } \right]}
\end{equation}
\noindent
and (\ref{even-odd}).
\end{proof}

\noindent
{\bf Note}. Corollary \ref{c-even} and the 
fact that the index $l$ is fixed yield the 
proof of Theorem \ref{odd-red}.

\medskip

From now on we assume $l \in \mathbb{N}$ is a fixed odd number. Consider now
the sets
\begin{equation}
V_{k,a}^{l} := 
\left\{ \nu_{2} \left( C_{l, 2^{k}(m-1)+a} \right): \, m \in \mathbb{N} \right\}
\end{equation}
\noindent
for $k \in \mathbb{N}$ and $1 \leq a \leq 2^{k}$.  Observe that, for 
fixed $k \in \mathbb{N}$ the $2^{k}$
sets $V_{k,a}^{l}$ contain all the information of the 
sequence $\{ \nu_{2} \left( C_{l,m} \right): \, m \geq 1 \}$. 
For example, for $k=2$, we have
\begin{eqnarray}
V_{2,1}^{l} &=& \left\{ \nu_{2}(C_{l, 4m-3}): m \in \mathbb{N}  \right\}, 
\nonumber \\
V_{2,2}^{l} &=& \left\{ \nu_{2}(C_{l, 4m-2}): m \in \mathbb{N}  \right\}, 
\nonumber \\
V_{2,3}^{l} &=& \left\{ \nu_{2}(C_{l, 4m-1}): m \in \mathbb{N}  \right\},
\nonumber \\
V_{2,4}^{l} &=& \left\{ \nu_{2}(C_{l, 4m}): m \in \mathbb{N}  \right\}.
\nonumber
\end{eqnarray}
\noindent
These four sets correspond to the third generation in the tree shown in
Figure \ref{tree5}. We also introduce the difference between $V_{k,a}^{l}$
and the basic sequence $\{ \nu_{2}(m): m \in \mathbb{N} \}$. 

\medskip

\noindent
{\bf Note}. The sets $V_{k,a}^{l}$ are attached to the vertices in the 
$k$-th generation of the decision tree $T(l)$. The terminal vertices of 
$T(l)$ are those corresponding
to indices $a$ such that the set 
\begin{equation}
S_{k,a}^{l} := 
\left\{ \nu_{2} \left( C_{l, 2^{k}(m-1)+a} \right) - \nu_{2}(m): \, m \in \mathbb{N} \right\}
\end{equation}
\noindent
reduces to a single value. \\

\begin{definition}
Let $k^{*}(l)$ be the first generation in the decision tree $T(l)$ that
contains a terminal vertex. This is the minimal $k$ for which there 
exist an index $a$, 
in the range $1 \leq a \leq 2^{k}$, such that
$V_{k,a}^{l}$ is a constant shift of the sequence $\{ \nu_{2}(n): n \in 
\mathbb{N} \}$. 
\end{definition}


\medskip

\begin{theorem} \label{thm_main}
Let $l \in \mathbb{N}$ be odd. Then $k^{*}(l) = \lfloor{ \log_{2}l \rfloor}$.
The $k^{*}-$th generation  contains $2^{k^{*}+1}-l$ terminal vertices. The tree 
$T(l)$ has one more generation consisting of $2(l-2^{k^{*}})$ terminal 
vertices. And these are the only terminal vertices.
\end{theorem}

The proof is divided into a sequence of steps.

\begin{lemma} \label{lem_1st}
Let $l$ be an odd integer and define $k$ via $2^{k} < l < 2^{k+1}$. Then
for $a$ in the range $1 \leq a \leq 2^{k+1}-l$ define 
$j_{1}(l,k,a):= -l+2(1+2^{k}-a)$. Then 
\begin{equation}
\nu_{2} \left( C_{l,2^{k}(m-1)+a} \right) = 
\nu_{2}(m) + \gamma_1(l,k,a) 
\end{equation}
\noindent
with 
\begin{equation}
\gamma_1(l,k,a) = l+k+1 + \nu_{2} \left( (j_{1}+l-1)! \times (l-j_{1})! 
\right). 
\end{equation}
Therefore, the vertex corresponding to the index $a$ is a terminal vertex
for the tree $T(l)$ with vertex constant $\gamma_1(l,k,a)$.
\end{lemma}
\begin{proof}
We have 
\begin{eqnarray}
\nu_{2} \left( C_{l,2^{k}(m-1)+a} \right) & = & 
\nu_{2} \left( A_{l,l+ 2(2^{k}(m-1)+a-1)} \right)  \label{eqn1} \\
& = & l + \sum_{j=-l+1}^{l} \nu_{2} \left( l + 2(2^{k}(m-1)+a-1) + j 
\right). \nonumber 
\end{eqnarray}
\noindent
The bounds on $a$ imply that $2-l \leq j_{1} \leq 2^{k+1}-l$ showing that
$j_{1}$ is in the range of summation. Morever it isolates the term
$2^{k+1}m$; that is, (\ref{eqn1}) can be computed as
\begin{eqnarray}
\nu_{2} \left( C_{l,2^{k}(m-1)+a} \right) & = & 
 l + \sum_{b=1}^{j_{1}+l-1} \nu_{2}(2^{k+1}m -b) + 
k+1 + \nu_{2} (m) + 
\sum_{b=1}^{l-j_{1}} \nu_{2}(2^{k+1}m + b ). \nonumber
\end{eqnarray}
\noindent
In the first sum we have 
$b \leq j_{1} + l-1 = 1-2a + 2^{k+1} < 2^{k+1}$, and in the second one 
$b \leq l - j_{1} = 2 \left( l-1+a-2^{k})  \right) < 2^{k+1}$,
by the choice of the upper bound on $a$. We conclude that 
\begin{eqnarray}
\nu_{2} \left( C_{l,2^{k}(m-1)+a} \right) & = & 
\nu_{2}(m) + l + k + 1 + 
\sum_{b=1}^{j_{1}+l-1} \nu_{2}(b) + 
\sum_{b=1}^{l- j_{1}} \nu_{2}(b). 
\nonumber
\end{eqnarray}
\noindent
This is the stated result. 
\end{proof}

\begin{lemma}
Let $k$ and $l$ be defined as above, then
\begin{eqnarray} \label{lbl_pm}
\nu_{2} \left( A_{l,2^{k+1}m+a} \right) = \nu_{2} \left( A_{l,2^{k+1}m-a-1} 
\right),\end{eqnarray}
for any $m \ge 1$ and $0 \le a < 2^{k+1} - l$.
\end{lemma}
\begin{proof}
Since $2^{k} < a+l < 2^{k+1}$, $0 \le a \le 2^{k}$, and $-2^{k+1} < a - l < 0$. Therefore,
\begin{eqnarray*}
\nu_{2} \left( A_{l,2^{k+1}m+a} \right) 
	&=& l + \sum_{j=-l+1}^{l} \nu_{2} \left( 2^{k+1}m+a+j \right)	\\
	&=& l + \sum_{j=-l}^{l-1} \nu_{2} \left( 2^{k+1}m-a-j \right)	\\
&=& \nu_{2} \left( A_{l,2^{k+1}m-a-1} \right).
\end{eqnarray*}
\end{proof}

Now since the sets 
$\{2^{k+1}m \pm a\, |\, 0 \le a < 2^{k+1} + 2^k -l, m \ge 1 \}$ 
and $\{2^{k+2}m \pm b\, |\, 2^{k+1}-l < b < l, m \ge 1 \}$ partition 
the set $\{a\, |\, a \ge l\}$, to prove the second half of 
Theorem~\ref{thm_main}, we only need to show the following.
\begin{lemma} \label{lem_2nd}
Let $k$ and $l$ be defined as above, then for $a$ in the 
range $2^{k+1}-l < a \leq 2^{k+1}$ define 
$j_{2}(l,k,a):= -l+2(1+2^{k+1}-a)$. Then 
\begin{equation}
\nu_{2} \left( C_{l,2^{k}(m-1)+a} \right) = 
\nu_{2}(m) + \gamma_2(l,k,a) 
\end{equation}
\noindent
with 
\begin{equation}
\gamma_2(l,k,a) = l+k+2 + \nu_{2} \left( (j_{2}+l-1)! \times (l-j_{2})! 
\right). 
\end{equation}
Therefore, the vertex corresponding to the index $a$ is a terminal vertex
for the tree $T(l)$ with vertex constant $\gamma_2(l,k,a)$.
\end{lemma}
\begin{proof}
The proof is the same as that of Lemma~\ref{lem_1st}, and thus omitted.
\end{proof}
\begin{example}
In the case $l=3$ we can take $k=1$. On the higher level, the restrictions 
on the parameter $a$ 
imply that must have $a=1$. A direct calculation shows that 
$j_{1}(3,1,1) = 1$ and $\gamma_1(3,1,1) = 7$. For the bottom two 
vertices, $a=2,4$; and we have $j_2(3,1,2)=3, \gamma_2(3,1,2)=9$; while 
$j_2(3,1,4)=-1, \gamma_2(3,1,4)=9$. 
This confirms the data on Figure~\ref{dec-3}. 
\end{example}

\begin{example}
For $l=5$, the theorem predicts three terminal vertices at the level $k=2$, 
corresponding to the values $a=1, 2, 3$. This confirms Figure \ref{tree5}
with terminal values given by $\gamma_1(5,2,1) = \gamma_1(5,2,3) = 13$ and 
$\gamma_1(5,2,2) =14$.  Similar results can be drawn for the level $k=3$. 
As before, the tree produces an explicit formula for the 
$2$-adic valuation of $C_{5,m}$. Indeed, define
\begin{equation}
f_{5}(m) = \begin{cases}
14  +  \nu_{2} \left( \frac{m+2}{4} \right)  & \text{ if } m \equiv 2 \bmod 4, \\
13  +  \nu_{2} \left( \frac{m+1}{4} \right)  & \text{ if } m \equiv 3 \bmod 4, \\
13  +  \nu_{2} \left( \frac{m+3}{4} \right)  & \text{ if } m \equiv 1 \bmod 4, \\
16  +  \nu_{2} \left( \frac{m}{8} \right)  & \text{ if } m \equiv 0 \bmod 8, \\
16  +  \nu_{2} \left( \frac{m+4}{8} \right)  & \text{ if } m \equiv 4 \bmod 8. 
\end{cases}
\end{equation}
\noindent
then, 
\begin{equation}
\nu_{2}(C_{5,m}) = f_{5}(m). 
\end{equation}
\end{example}
To finish the proof of Theorem~\ref{thm_main}, we need to establish

\begin{lemma}
There are no terminal vertices of level less than $k$.
\end{lemma}
\begin{proof}
The value of a vertex on the level $u < k$ is obtained from
$\nu_{2} \left( C_{l,2^{u}(m-1)+a} \right)$. The proof of 
Lemma~\ref{lem_1st}, shows that
\begin{equation}
\nu_{2} \left( C_{l,2^{u}(m-1)+a} \right) = \sum_{i=0}^{v}\nu_{2}(m+i) + c,
\nonumber
\end{equation}
\noindent
for some constants $v > 0$ and $c$. The next lemma proves that this cannot 
happen.
\end{proof}

\begin{lemma}
If 
\begin{equation}
\sum_{i=0}^{a} \nu_{2}(m+i) = \nu_{2}(m+b) + c
\end{equation}
\noindent
for all $m \ge 1$, 
and some constants $a, b, c$, then $a = b = c = 0$.
\end{lemma}
\begin{proof}
Suppose the lemma is not true and $b>a$. 
Choose $m$ such that $m+b = 2^u$ for some $u$, then
$\sum_{i=0}^{a} \nu_{2}(m+i) = \nu_{2}((b-a) \cdots b)$. Therefore
$c = \nu_{2}((b-a) \cdots b) - u$. Similarly choose $m$ such that 
$m+b = 2^{u+1}$, and conclude that 
$c = \nu_{2}((b-a) \cdots b) - u - 1$. This is a contradiction. 
The proof to the other two cases where $0 \le b \le a$ and $b < 0$ are 
similar, and thus omitted.
\end{proof}

\medskip

The final figure shows how to produce the trees corresponding to $l$ odd. 
First determine $n$ by $2^{n} < l < 2^{n+1}$ and form a complete binary 
tree $T$ where the last level has $2^{n}$ vertices. Now from $T$ branch an
odd number of vertices that yields the decision trees $T(l)$. 
Figure \ref{dec-odda} shows the four trees 
corresponding to the odd indices $l$ in the range $8 < l < 16$.  \\

{{
\begin{figure}[ht]
\begin{center}
\centerline{
\psfig{file=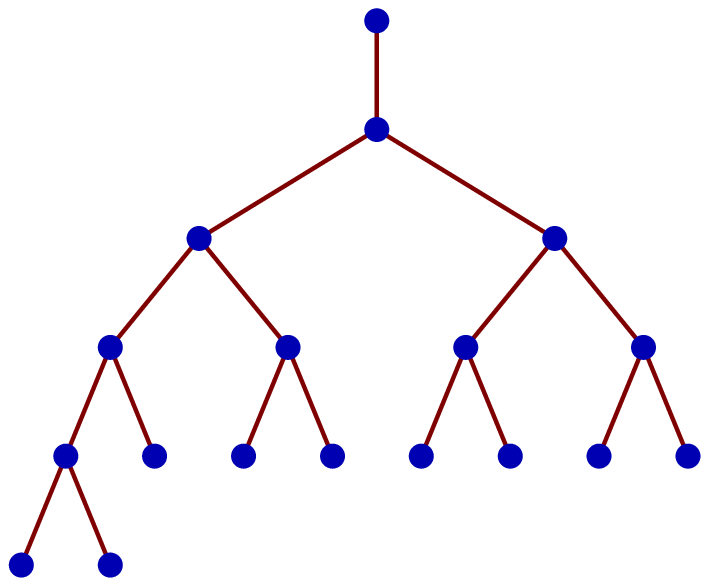,width=15em,angle=0}
\psfig{file=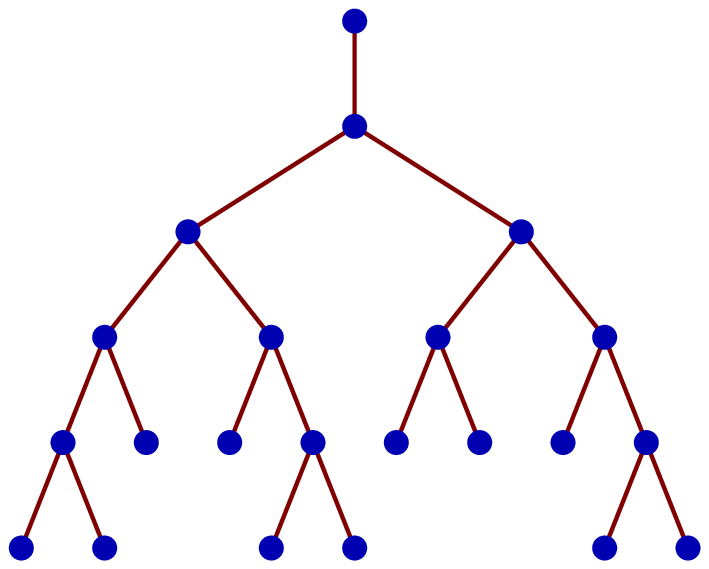,width=15em,angle=0}
}
\label{dec-odda1}
\end{center}
\end{figure}
}}

{{
\begin{figure}[ht]
\begin{center}
\centerline{
\psfig{file=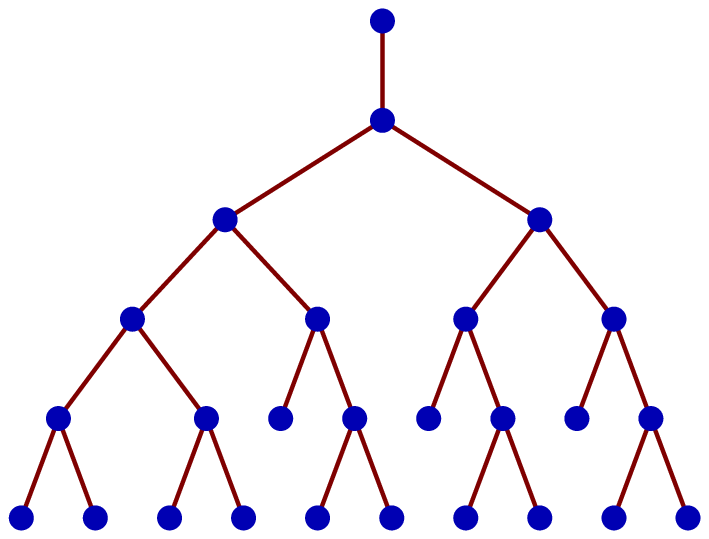,width=15em,angle=0}
\psfig{file=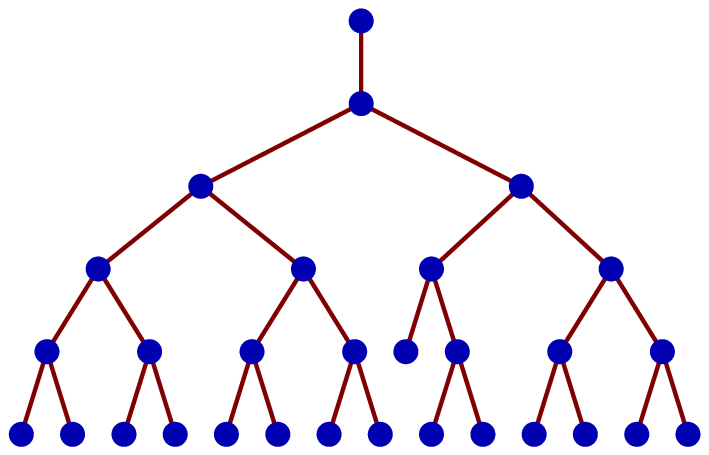,width=15em,angle=0}
}
\caption{The trees for $l$ odd between $8$ and $16$}
\label{dec-odda}
\end{center}
\end{figure}
}}

\noindent
{\bf Acknowledgements}. The second author acknowledges the partial support of 
NSF-DMS 0713836.

\end{document}